\documentclass{article}
\usepackage{amssymb,amsmath,graphicx,ucs,hyperref}
\usepackage[utf8x]{inputenc}

\newtheorem{theorem}{Theorem}

\title{Compact packings of the space\\ with two sizes of spheres}
\author{Thomas Fernique\footnote{Univ. Paris 13, CNRS, Sorbonne Paris Cit\'e, UMR 7030, 93430 Villetaneuse, France.}}
\date{}

\begin{document}

\maketitle

\begin{abstract}
Compact packings are specific packings of spheres which can be seen as tilings and are good candidates to maximize the density.
We show that the compact packings of the Euclidean space with two sizes of spheres are exactly those obtained by filling with spheres of size $\sqrt{2}-1$ the octahedral holes of a close-packing of spheres of size $1$.
\end{abstract}

\section{Introduction}

A {\em sphere packing} is a set of interior-disjoint spheres.
Its {\em density} is the upper limit for $r\to \infty$ of the proportion of the space within distance $r$ from the origin which is covered by the interiors of the spheres.
We shall here always consider packings up to similarity (the density is unchanged).
Finding the densest packings of spheres of given sizes is an issue of great interest in chemistry in order to model or predict new material structures (see, {\em e.g.}, \cite{HO11}).\\

In \cite{FT64}, a packing of the Euclidean plane by circles is said to be {\em compact} if its {\em contact graph}, {\em i.e.}, the graph which connects the centers of adjacent circles, is a triangulation.
There is only one compact packing of the plane with equal circles, namely the {\em hexagonal packing} (circles are centered on a triangular grid).
It is known to maximize the density among packing by equal circles \cite{Thu10,FT43}.
There are $9$ pairs $(1,r)$ allowing compact packings of the plane by circles of size $1$ and $r$ \cite{Ken06} and $164$ triples $(1,r,s)$ allowing compact packings of the plane by circles of size $1$, $r$ and $s$ \cite{FHS18}.
For $7$ of the $9$ ratio allowing a compact packing by two circles, the maximal density has been proved to be reached by a compact packing \cite{Hep00,Hep03,Ken04}.
The two other cases are still open.\\

By analogy, a packing of spheres in $\mathbb{R}^n$ is said to be {\em compact} if its contact graph is a homogeneous simplicial complex of dimension $n$, {\em i.e.}, the $1$-skeleton of a covering of $\mathbb{R}^n$ by interior-disjoint $n$-dimensional simplices, the intersection of any two simplices being either empty or a lower dimension face (such a covering is also called a {\em face-to-face tiling}).\\

Compact packings by equal spheres are known to exist in dimension $8$ and $24$ (spheres are respectively centered on the $E_8$ and Leech lattices).
In both cases, they maximize the density among packings by equal spheres \cite{Via17,CKMRV17}.
In dimension 3, however, there is no compact packing by equal spheres.
It is indeed impossible to tile the Euclidean space by regular tetrahedra.
What about packings with two sizes of spheres?
Are there compact ones?
If so, can we find them all?
This is the issue addressed here.\\

\begin{figure}[hbtp]
\centering
\includegraphics[width=0.95\textwidth]{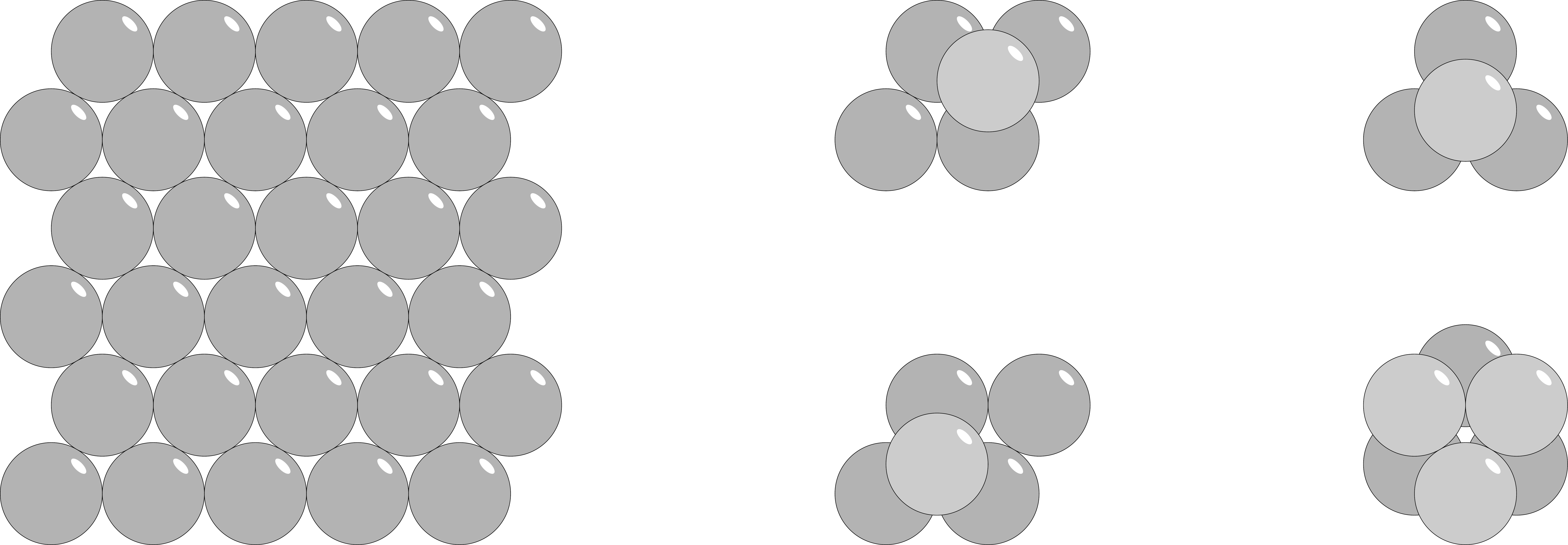}
\caption{
Layer of a close-packing (left).
Only one out of two neighbor voids between three spheres of a layer can be filled by a sphere of the next layer (center, top view).
Close-packing have two types of holes: tetrahedral around a filled void (top-right) and octahedral around a non-filled void (bottom-right).
}
\label{fig:close_packing}
\end{figure}

The maximal density of packing by equal spheres in the Euclidean space is achieved for the so-called {\em close-packings of equal spheres}, as conjectured by Kepler in 1611 and proven nearly four centuries later \cite{Hal05}.
Close-packings of spheres are formed by stacked layers of spheres centered on a triangular grid, where each sphere of a layer fills the void between three mutually adjacent spheres of the previous layer (Fig.~\ref{fig:close_packing}).
Only one void out of two of a layer can be filled by the spheres of the next one, and there is exactly two ways to put a layer on the last layer: either lined up with the last but one layer or not.
This yields uncountably many different packings, all with the same density.
The packing where each layer is {\em never} (resp. {\em always}) lined up with the last but one is said to be {\em cubic close-packed} or {\em face-centered cubic} (resp. {\em hexagonal close-packed}).\\

In a packing, a {\em hole} is a local maxima of the distance to the closest sphere.
Compact packings have only {\em simplicial holes}, {\em i.e.}, holes at equal distance from $n+1$ spheres in $\mathbb{R}^n$.
Close-packings of unit spheres in $\mathbb{R}^3$ do have simplicial (tetrahedral) holes between three adjacent spheres in a layer and the sphere which fills the void between them.
But half of these voids are not filled and create so-called {\em octahedral holes} between three spheres of a layer and three spheres of the next layer (Fig.~\ref{fig:close_packing}).
Inserting small spheres centered on these octahedral holes turns out to be the one and only way to get compact packings:

\begin{theorem}\label{th:main}
The compact packings of the Euclidean space with two sizes of spheres are exactly those obtained by filling with spheres of size $\sqrt{2}-1$ the octahedral holes of a close-packing of spheres of size $1$.
\end{theorem}

Such structures are well known in chemistry.
In particular, filling the octahedral holes of a face-centered cubic packing yields the rock salt structure.
The main point here is to show that there is no other compact packing.
This may be a bit disappointing in comparison with the variety of compact packings by two or three discs in the Euclidean plane, but the method opens the way for three sizes of spheres, where a greater variety may be expected (see, {\em e.g.}, the discussion in \cite{LH93}).\\

These compact packings improve by a factor $5\sqrt{2}-6\simeq 1.07$ the density of close-packings by equal spheres.
The various results mentioned above suggest that compact packings, when they exist, are good to maximize the density among packings by the same spheres.
In this case, this issue is still open.\\

The rest of this paper is organized as follows.
Section~\ref{sec:necklaces} introduces a first type of local configuration called {\em necklace}.
We show in Section~\ref{sec:skew_necklaces} how so-called {\em skew} necklaces allow to reduce to $10$ the number of possible values for the ratio of sphere sizes.
We use other necklaces in Section~\ref{sec:large_necklaces} to reduce this number to one, namely $\sqrt{2}-1$.
The double inclusion of Theorem~\ref{th:main} is proven in a combinatorial way in Sections~\ref{sec:proof1} and \ref{sec:proof2}, where a second type of local configuration called {\em shell} is introduced.



\section{Necklaces}
\label{sec:necklaces}

Any pair of adjacent spheres $B$ and $H$ in a compact packing yields a sequence $B_1,\ldots,B_k$ of spheres, where $B_k$, $B_{k+1}$, $B$ and $H$ are pairwise tangent.
We call such a sequence a {\em necklace}, seing $B$ as the body, $H$ as the head and the $B_k$'s as the beads.
A necklace is coded by the word over $\{1,r\}$ whose $k$-th letter gives the radius of the $k$-th bead.
Among the possible codings, we usually choose the lexicographically minimal one.
A necklace is said to be {\em large} if $B$ and $H$ have radius $1$, {\em small} if they have radius $r$, and {\em skew} otherwise.
Any compact packing by spheres of radius $1$ and $r$ contains two different spheres in contact, thus at least one skew necklace.

\section{Skew necklaces and possible values of $r$}
\label{sec:skew_necklaces}

Any skew necklace contains at most $5$ beads.
Indeed, the spheres of a necklace can be rolled on $H$ until their center are all in the same plane as the center of $H$.
Since, in the plane, at most $6$ unit discs can be adjacent to a unit discs, there is at most $6$ beads, actually even at most $5$ beads because rolling the beads on $H$ enlarged the necklace.
There is thus finitely many potential skew necklaces, see Tab.~\ref{tab:skew_necklaces}.
We shall see that each skew necklace characterizes the radius $r$.\\

\begin{table}[hbtp]
\centering
\begin{tabular}{llllll}
11111 & 1111r & 111rr & 11rrr & 1rrrr & rrrrr\\
1111 & 111r & 11r1r & 1r1rr & rrrr &\\
111 & 11r & 11rr & 1rrr & &\\
&& 1r1r & rrr & &\\
&& 1rr & & &
\end{tabular}
\caption{Potential skew necklaces (same number of r in each column).}
\label{tab:skew_necklaces}
\end{table}

In a necklace, consider the tetrahedron whose vertices are the centers of $B$, $H$, $B_{k}$ and $B_{k+1}$.
Denote by $\delta_{s,t}$ the dihedral angle between the faces $BHB_{k}$ and $BHB_{k+1}$, where $s$ and $t$ are the radius of $B_k$ and $B_{k+1}$.
Let $\widehat{xyz}$ denote the angle in $y$ of a contact graph between three spheres of radius $x$, $y$ and $z$.
Spherical geometry allows to express the cosine of $\delta_{s,t}$ as a function of the angles of faces of the tetrahedron:
$$
\cos \delta_{s,t}=\frac{\cos\widehat{s1t}-\cos\widehat{r1s}\cos\widehat{r1t}}{\sin\widehat{r1s}\sin\widehat{r1t}}.
$$
Sines can be expressed with cosines via $\sin^2a+\cos^2a=1$, and cosines can be expressed with $r$ via the cosine theorem:
$$
\cos\widehat{111}=\frac{1}{2},\qquad
\cos\widehat{11r}=\frac{1}{1+r},\qquad
\cos\widehat{r1r}=1-\frac{2r^2}{(1+r)^2}.
$$
A computation then yields:
$$
\cos \delta_{1,1}=\frac{r^2+2r-1}{2r(2+r)},\quad
\cos \delta_{1,r}=\sqrt{\frac{r}{(2+r)(2r+1)}},\quad
\cos \delta_{r,r}=\frac{1+2r-r^2}{2(2r+1)}.
$$

The sum of the dihedral angles given by all the consecutive beads in a necklace must be equal to $2\pi$: this yields an equation in $r$.
For example, the skew necklace 111rr yields
$$
2\delta_{1,1}+2\delta_{1,r}+\delta_{r,r}=2\pi.
$$
Take the cosine of both sides of this equation, fully expand the left-hand side and substract $\cos(2\pi)=1$ to both sides.
This yields a polynomial equation in cosines and sines of the $\delta_{s,t}$'s.
Use $\sin^2a+\cos^2a=1$ to replace the sines by cosines, then the above expressions to get an equation in $r$.
This equation may contain square roots.
To avoid this, we introduce the auxiliary variables $X_0,\ldots, X_3$ defined by the algebraic equations below (left) and use them to express cosines and sines (right):\\
\begin{tabular}{c|c}
\begin{minipage}{0.49\textwidth}
\begin{eqnarray*}
0 &=& (2+r)(2r+1)X_0^2-r,\\
0 &=& X_1^2-3r^2-6r+1,\\
0 &=& (2+r)(1+2r)X_2^2-2,\\
0 &=& X_3^2+r^2-6r-3.\\
\end{eqnarray*}
\end{minipage}
&
\begin{minipage}{0.49\textwidth}
\begin{eqnarray*}
\cos \delta_{1,r}&=&X_0,\\
\sin \delta_{1,1}&=&\tfrac{1+r}{2r(2+r)}X_1,\\
\sin \delta_{1,r}&=&(1+r)X_2,\\
\sin \delta_{r,r}&=&-\tfrac{1+r}{4r+2}X_3.\\
\end{eqnarray*}
\end{minipage}
\end{tabular}
~\\
Thus, for each skew necklace we get a system of $5$ algebraic equations (one for the sum of dihedral angles and $4$ for the auxiliary variables) that we can solve exactly with a computer algebra software (we used \cite{sagemath}).
This yields a total of $16$ possible values of $r$.
For each one, we check by interval arithmetic that not only the cosine of the equation on dihedral angles is satisfied, but the equation itself: this amounts to check that the sum of the angles is equal to $2\pi$ and no $2k\pi$ for some $k\neq 1$.
This reduces to $10$ the number of possible values of $r$ (Tab.~\ref{tab:radii}).\\

\begin{table}[hbtp]
\centering
\begin{tabular}{lll}
11111 & $X^{4} + 4 X^{3} + X^{2} - 6 X + 1$ & $0.902$\\
1111r & $4 X^{4} + 8 X^{3} - 4 X^{2} - 6 X + 1$ & $0.849$\\
111rr & $X^{4} + 4 X^{3} + 3 X^{2} - 6 X + 1$ & $0.720$\\
11r1r & $4 X^{3} - 20 X^{2} + 9 X + 2$ & $0.690$\\
11rrr & $X^{4} - 2 X^{3} - 5 X^{2} + 1$ & $0.420$\\
1111  & $X^{2} + 2 X - 1$ & $0.414$\\
111r  & $2 X^{2} + 3 X - 1$ & $0.280$\\
111   & $2 X^{2} + 4 X - 1$ & $0.224$\\
1r1rr & $2 X^{3} + 9 X^{2} - 20 X + 4$ & $0.223$\\
11rr  & $X^{2} - 6 X + 1$ & $0.171$
\end{tabular}
\caption{Skew necklace (left), minimal polynomial of the value of $r$ it characterizes (middle) and approximated value of $r$ (right).}
\label{tab:radii}
\end{table}

\section{Large and small necklaces}
\label{sec:large_necklaces}

Let $\overline{\delta}_{s,t}$ be the analogous of $\delta_{s,t}$ for large necklaces, {\em i.e.}, when $B$ and $H$ both have radius $1$.
A computation yields (as for $\cos\delta_{s,t}$):
$$
\cos\overline{\delta}_{1,1}=\frac{1}{3},\qquad
\cos\overline{\delta}_{1,r}=\frac{1}{\sqrt{3r(2+r)}},\qquad
\cos\overline{\delta}_{r,r}=\frac{2-r}{2+r}.
$$
The large necklaces correspond to the non-zero triples $(i,j,k)$ of non-negative integers satisfying
$$
i\overline{\delta}_{1,1}+j\overline{\delta}_{1,r}+k\overline{\delta}_{r,r}=2\pi.
$$
Given a value $r$, we first compute a lower approximation of the $\overline{\delta}_{s,t}$'s which yields an upper bound on the possible values for $i$, $j$ and $k$.
We then use arithmetic interval to find triples $(i,j,k)$ which could be solutions.
Last, we check these triples exactly by taking the cosine of the equation, expanding it and use the exact values of the $\cos\overline{\delta}_{s,t}$'s (and the sines via $\sin^2a+\cos^2a=1$).
We find that only $r=\sqrt{2}-1$ allows large necklaces, namely 111r1r and 11r11r.\\

Last, let $\underline{\delta}_{s,t}$ be the analogous of $\delta_{s,t}$ for small necklaces, {\em i.e.}, when $B$ and $H$ both have radius $r$.
We compute\footnote{$\cos\underline{\delta}_{1,1}$ and $\cos\underline{\delta}_{1,r}$ can be derived from $\cos\overline{\delta}_{r,r}$ and $\cos\overline{\delta}_{1,r}$ by replacing $r$ by $\tfrac{1}{r}$: it amounts to consider the same dihedral angle in homothetic tetrahedra.}:
$$
\cos\underline{\delta}_{1,1}=\frac{2r-1}{2r+1},\qquad
\cos\underline{\delta}_{1,r}=\frac{r}{\sqrt{3(1+2r)}},\qquad
\cos\underline{\delta}_{r,r}=\frac{1}{3}.
$$
Then we search for non-zero triples $(i,j,k)$ of non-negative integers satisfying
$$
i\underline{\delta}_{1,1}+j\underline{\delta}_{1,r}+k\underline{\delta}_{r,r}=2\pi.
$$
We find that only $r=3-2\sqrt{2}$ allows small necklaces, namely 11rr.

\section{First inclusion}
\label{sec:proof1}

Consider a compact packing by spheres of size $1$ and $r\in(0,1)$.
It must contains a large sphere in contact with a small one, hence a skew necklace.
There are only $10$ possible skew necklaces.
Since any of them contains two adjacent beads of radius $1$ (see Tab.~\ref{tab:radii}), there must be two adjacent spheres of radius $1$ in the packing, thus a large necklace.
The only $r$ which allows a large necklace is $r=\sqrt{2}-1$.
This is thus the only one which may allow compact packings.\\

For $r=\sqrt{2}-1$, there is no small necklace and only one skew necklace: 1111.
A small sphere can thus be surrounded only by large spheres.
Consider the tetrahedra which connect the center of such a small sphere to the centers of three mutually adjacent large spheres.
The solid angle in the small sphere of these tetrahedra can be computed from $r$ using the Girard's, cosine and spherical cosine theorems.
For $r=\sqrt{2}-1$, this yields $\tfrac{\pi}{2}$.
There are thus $8$ such tetrahedra around each small sphere, {\em i.e.}, $6$ large spheres centered on the vertices of a regular octahedra.\\

Spheres of size $\sqrt{2}-1$ can thus be inserted in the octahedral holes of any close-packing by spheres of size $1$.
This transforms all the octahedral holes into $8$ tetrahedral holes, {\em i.e.}, it yields a compact packing.
This proves the first inclusion of Theorem~\ref{th:main}: filling the octahedral holes of a close-packing yields a compact packing.

\section{Second inclusion}
\label{sec:proof2}

We here prove the second inclusion of Theorem~\ref{th:main}: any compact packing can be obtained by filling the octahedral holes of some close-packing.
Fix a compact packing by spheres of size $1$ and $\sqrt{2}-1$.
Consider a large sphere $S$ in this packing.
Let us call {\em shell} the set of spheres which are adjacent to it and represent it by a spherical triangulation with vertices labelled in $\{1,r,s\}$: the vertices represent the spheres and the edges their contact graph.
This shell cannot contains only large spheres, because the solid angle of a regular tetrahedron does not divide $4\pi$.
It thus contains a small sphere.
In the shell, this small sphere is surrounded by four large spheres because the unique skew necklace is 1111.
We make two cases.\\

\begin{figure}[hbtp]
\centering
\includegraphics[width=\textwidth]{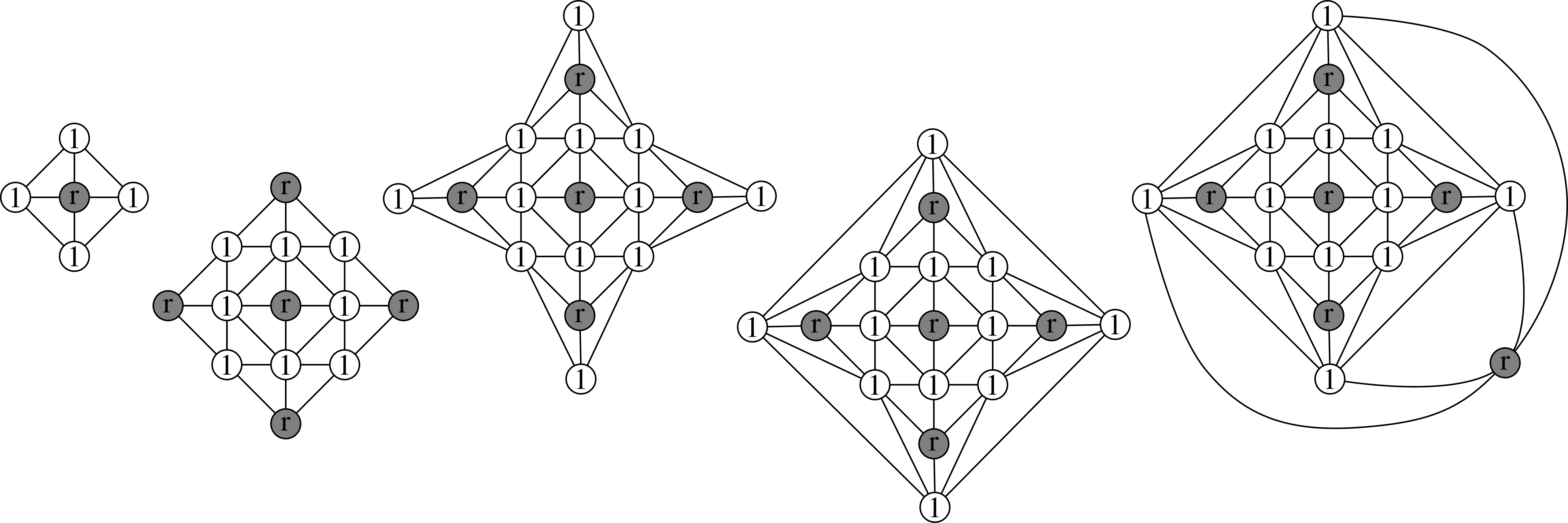}
\caption{Finding the shell of a large sphere: first case.}
\label{fig:large_shell_1}
\end{figure}

{\bf First case}.
Assume that each of these four large spheres is surrounded by 11r11r.
We proceed in $5$ steps illustrated on Fig.~\ref{fig:large_shell_1}.
\begin{enumerate}
\item Start from a small sphere surrounded by four large ones.
\item By hypothesis, each large sphere is surrounded by 11r11r.
\item Each new small sphere is surrounded by four large spheres.
There are now $12$ large spheres.
\item There are four large spheres already surrounded by 11r11r.
Since the unique large necklaces are 11r11r and 111r1r, they have no other neighbor.
\item Since no thirteen large sphere can be added \cite{SW53}, there must be a new small sphere which has to be connected with the four large spheres whose neighborhoods were still uncompleted.
Each sphere has a complete neighborhood.
The spherical triangulation is completed.
\end{enumerate}

\begin{figure}[hbtp]
\centering
\includegraphics[width=\textwidth]{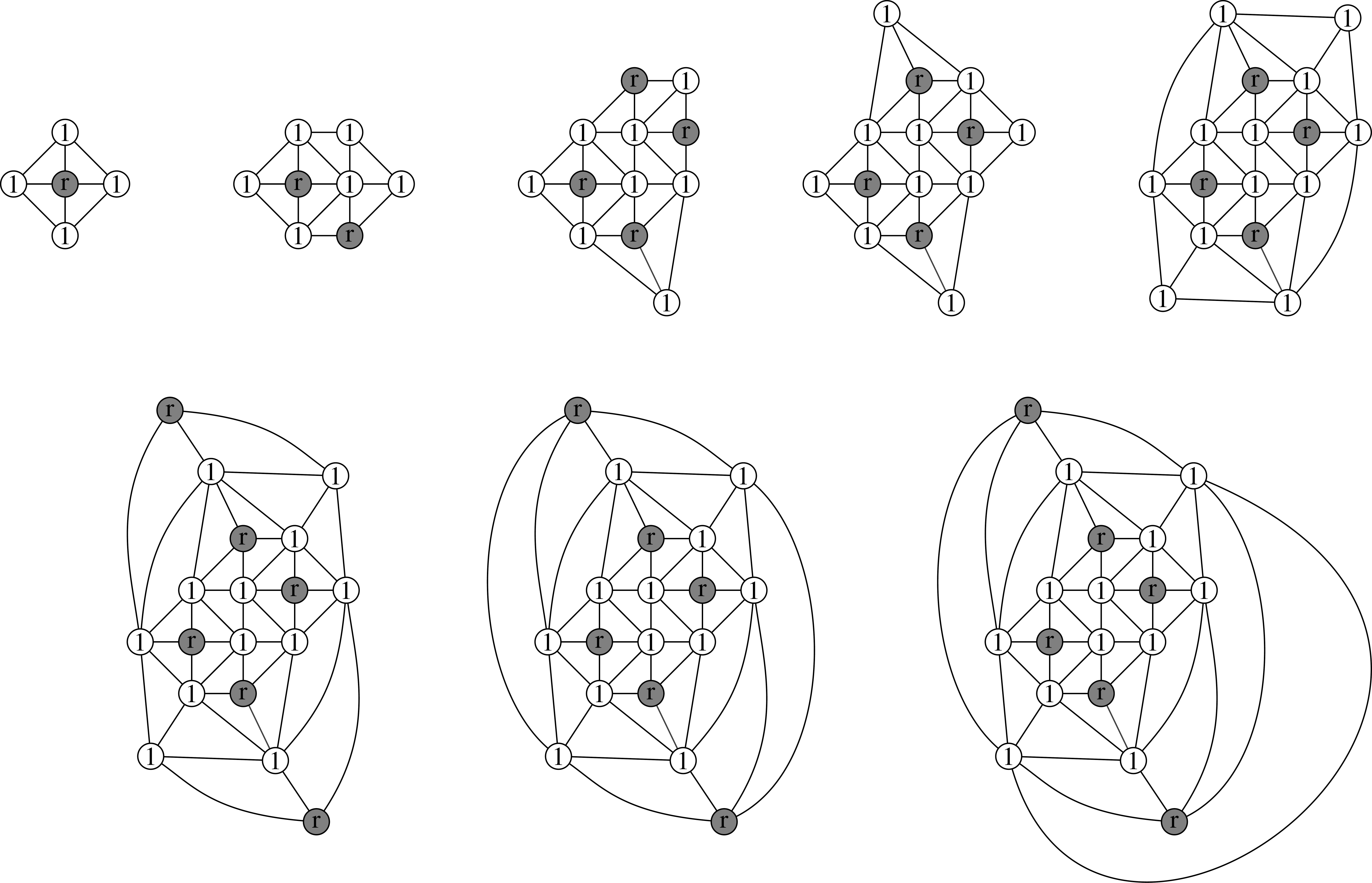}
\caption{Finding the shell of a large sphere: second case.}
\label{fig:large_shell_2}
\end{figure}

{\bf Second case}.
Assume that at least one of the four large spheres is not surrounded by 11r11r.
It is thus surrounded by 111r1r because these are the two unique large necklaces.
We proceed in $8$ steps illustrated on Fig.~\ref{fig:large_shell_2}.
\begin{enumerate}
\item Start from a small sphere surrounded by four large ones.
\item By hypothesis, there is a large sphere surrounded by 111r1r.
\item The small sphere is surrounded by 1111 (unique skew necklace).
The large sphere whose neighborhood contains three consecutive large spheres has to be surrounded by 111r1r, the unique compatible large necklace.
\item The two new small spheres must be surrounded by 1111.
\item Two large spheres are surrounded by 1r11r1, hence have a completed neighborhood.
Two large spheres are surrounded by 1r1r1, hence must be surrounded by 111r1r: this adds two large spheres.
There are now $12$ large spheres.
\item The two large spheres surrounded by 11r11 must be surrounded by 11r11r.
This adds two small spheres.
\item Two large spheres are surrounded by 11r11r, hence have a completed neighborhood..
\item Only two large spheres have an uncompleted neighborhood: they are surrounded by r111r.
They must thus be surrounded by 111r1r, {\em i.e.}, they must have a new large sphere as a neighbor.
Since no thirteen large sphere can be added \cite{SW53}, the only possibility is that these two large spheres are mutual neighbors.
Each sphere has a complete neighborhood.
The spherical triangulation is completed.
\end{enumerate}

\begin{figure}[hbtp]
\centering
\includegraphics[width=0.8\textwidth]{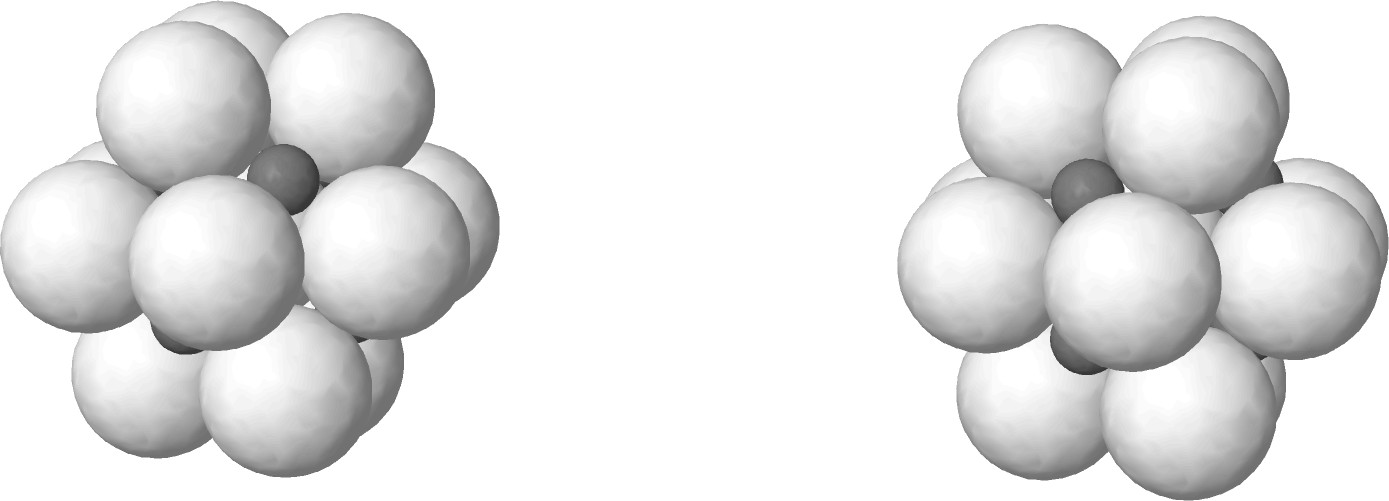}
\caption{
The two ways a large sphere can be surrounded in a compact packing by spheres of size $1$ and $\sqrt{2}-1$.
In each case there are $12$ large spheres and $6$ small ones.
In one case, the large spheres are centered on the vertices of a cuboctahedron and the small ones on the vertices of an octahedron (left).
In the other case, the large spheres are centered on the vertices of a triangular orthobicupola and the small ones on the vertices of a triangular prism (right).
}
\label{fig:large_shells}
\end{figure}

The two possible shells surrounding a large sphere are depicted on Fig.~\ref{fig:large_shells}.
Both have the following property: whenever the central large sphere has three consecutive neighbors which are coplanar large spheres, then it is surrounded by a ring of six coplanar large spheres (in the first case there are two such rings, in the second one only one).\\

\begin{figure}[hbtp]
\centering
\includegraphics[width=0.8\textwidth]{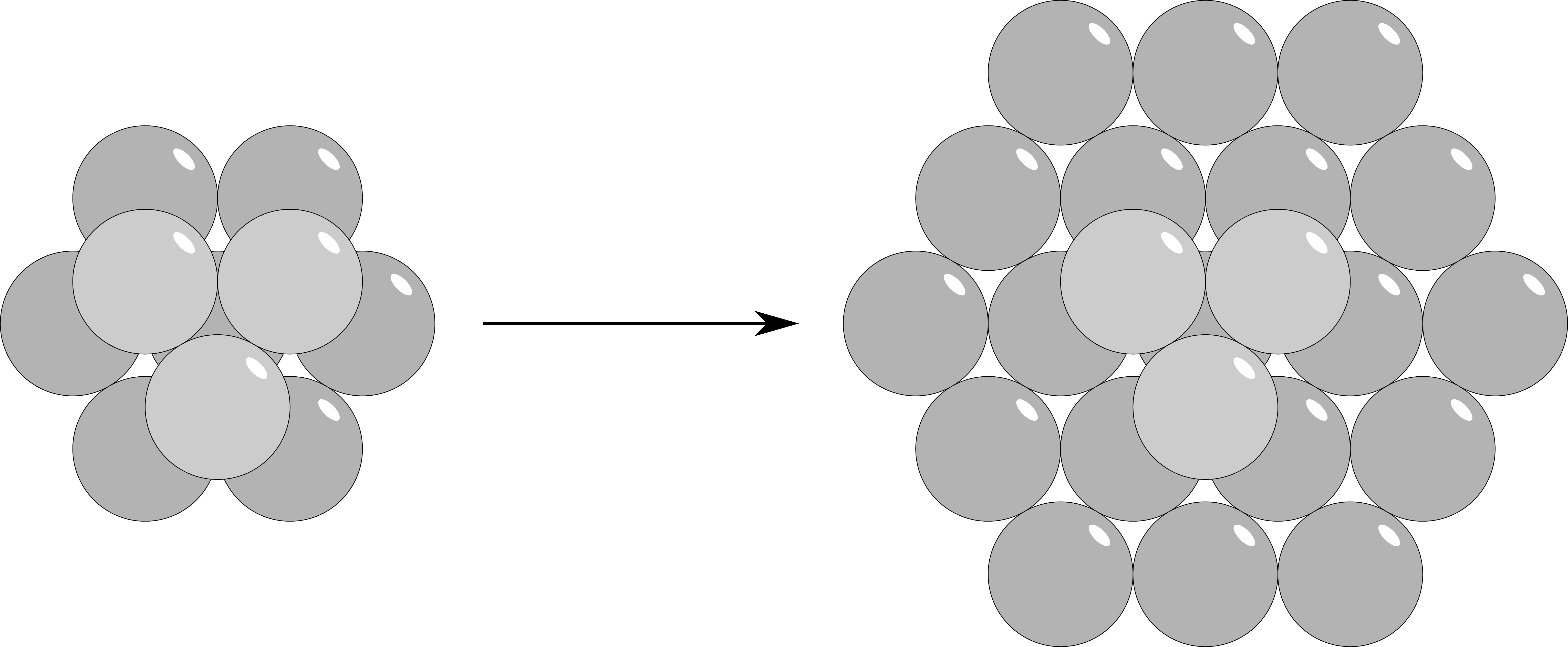}
\caption{
A large sphere, its shell and a ring of six coplanar spheres (left, top view, with coplanar spheres similarly shaded).
Each sphere of the ring has itself a ring (right).
This can be iterated to grow a whole layer of a close-packing.
}
\label{fig:layer_growth}
\end{figure}

Now, consider in a shell a ring of six coplanar large spheres around the central sphere (Fig.~\ref{fig:layer_growth}, left).
Each of these six spheres itself has three consecutive neighbors which are coplanar large spheres (two in the ring and the central sphere), so that they also are surrounded by a ring of six coplanar large sphere (Fig.~\ref{fig:layer_growth}, right).
Iterating this argument show that the initial sphere belongs to a set of coplanar large spheres centered on the vertices of a triangular grid, {\em i.e.}, a layer of a close-packing. 
Moreover, all the large spheres in this layer must have the same shell: each type of shell indeed fixes the shift between consecutive layers (in the first case the next layer is not line up with the previous one, in the second case both the next and previous layer are line up).
The large spheres of the considered compact packing thus form a close-packing.
Filling the octahedral holes of this close-packing by small spheres yields the compact packing itself.
This proves the second inclusion of Theorem~\ref{th:main}.

\bibliographystyle{alpha}
\bibliography{2balls}

\end{document}